\documentclass[10pt]{amsart}

\usepackage{amsfonts}
\usepackage{pdfsync}
\usepackage{graphicx}
\usepackage{amsmath}
\usepackage{setspace}
\usepackage[left=3cm,right=3cm,top=2.5cm,bottom=3cm]{geometry}
\usepackage{amsmath, %amssymb,
amsthm,bbm,color,graphics,version}
\usepackage{mathrsfs}

\usepackage[cp1250]{inputenc}
\usepackage{polski}

\theoremstyle{plain}
\theoremstyle{plain}

\def\m{\mu \;} %(the drift, to be decided ...)
\def\mbb{\mathbb}  
\def\mc{\mathcal} \def\unl{\underline} 

\def\P{{\mathbb P}} \def\le{\left} \def\ri{\right} \def\i{\infty}
\def\E{{\mathbb E}}

\def\te#1{\mathrm{e}^{#1}}

\def\I{\int}      
  \def\d{\delta}   \def\th{\theta}

     \def\s{\sigma}
     \def\ps{\psi}
 \def\q{\qquad} 
\def\F{\Phi}   
  \def\td{\text{\rm d}}

\numberwithin{equation}{section}

\newtheorem{Thm}{Twierdzenie}[section]
\newtheorem{Cor}[Thm]{Wniosek} 
 
{%\theorembodyfont{\normalfont}
 \newtheorem{Def}[Thm]{Definnicja}
\newtheorem{Rem}[Thm]{Uwaga} 
 
}

\title{Problem optymalizacyjny de Finettiego\\ dla procesów L\'evy'ego}

\author{Zbigniew Palmowski}
\address{Instytut Matematyczny, Uniwersytet Wrocławski.}
\email{zbigniew.palmowski@uwr.edu.pl}

\date{\today}

\begin{document}
\begin{abstract}
Celem tej pracy jest przedstawienie wybranych metod probabilistycznych w rozwiązywaniu tzw. dywidendowego problemu optymalizacyjnego
oraz pokazanie związków z innymi metodami analitycznymi związanymi m.in. z teorią równań różniczkowych.
\vspace{3mm}

\noindent {\sc Słowa kluczowe.} problem dywidendowy $\star$ proces L\'evy'ego $\star$ optymalna strategia $\star$ równanie HJB

\end{abstract}

\maketitle

\pagestyle{myheadings} \markboth{\sc Z.\ Palmowski}
{\sc Problem optymalizacyjny de Finettiego}

\vspace{1.8cm}

\tableofcontents

\newpage

\section{Przedstawienie problemu}\label{sec:model}
Teoria optymalizacji (dyskretnej i losowej) czy też teoria sterowania  stanowi obszerną dziedzinę matematyki związaną z równaniami
różniczkowymi, analizą matematyczną czy też informatyką. Celem tej pracy jest opis bardzo szczególnego stochastycznego problemu
optymalizacyjnego, który w ostatnim czasie skupił ogromną uwagę matematyków aktuarialno-finansowych oraz pokazanie jakie metody badawcze
i pytania są z nim związane. To na co pragniemy zwrócić szczególną uwagę
to są nowe metody probabilistyczne związane z martyngałami pozwalające na analizę procesów posiadających skoki.

W pracy będziemy zajmować się spektralnie ujemnym procesem L\'evy’ego $\{X_t,\, t\geq 0\}$. Jest to proces stochastyczny
o wartościach rzeczywistych, o stacjonarnych i niezależnych przyrostach
i o trajektoriach prawie wszędzie typu c\'adl\'ag\footnote{c\'adl\'ag= prawostronnie ciągły ze skończonymi lewostronnymi granicami}.
Określenie {\it spektralnie ujemny} oznacza, że dopuszczamy tylko skoki ujemne lub brak skoków.
Będziemy dopuszczać też możliwość dowolnej początkowej wartości $X_0=x$, co będzie podkreślone w oznaczeniu przestrzeni probabilistycznej, na której formalnie ten
proces jest zdefiniowany $(\Omega, \mathcal{F}, \{\mathcal{F}_t\}_{t\geq 0\}},\P_x)$, gdzie $\mathcal{F}_t$ jest naturalną, prawostronnie ciągłą filtracją procesu $X$
oraz $\mathcal{F}=\bigvee_{\{t\geq 0\}}\mathcal{F}_t$.
Później będziemy pisać $\P_0=\P$. Dwa przykłady procesów są dla nas najistotniejsze. Pierwszy to {\it liniowy ruch Browna}
\begin{equation}\label{Brown}
X_t=x+\eta t+\sigma B_t,\qquad \eta\in \mathbb{R},\sigma >0,
\end{equation}
gdzie $B_t$ jest ruchem Browna. Drugim przykładem jest proces Cram\'era-Lundberga
\begin{equation}\label{cramermod} X_t=x+ \eta t- \sum_{k=1}^{N_t} (C_k -\E C_1),
\end{equation}
gdzie $C_k$ tworzą ciąg nieujemnych, niezależnych zmiennych losowych o jednakowym rozkładzie
przychodzących zgodnie z niezależnym od nich, jednorodnym procesem Poissona z intensywnością $\lambda$.
Powyżej $\E_x[\cdot]$ ($\E[\cdot]$) jest wartością oczekiwaną względem $\P_x$ ($\P$) (często też będziemy pisać $\E_x[\cdot, \mathbf{1}_{A}]=\E_x[\cdot; A]$).
W teorii ryzyka oba procesy opisują ilość rezerw jakie firma posiada w momencie $t$.
W przypadku (\ref{cramermod}) $x$ jest początkowym kapitałem, zmienne $C_k$ opisują zgłaszane szkody, zaś $\eta+\lambda \E C_1$ podaje stałą intensywność wpłat.
Zwykle zakłada się, że $\E X_1>0$ (lub powyżej $\eta>0$) czyli, że rezerwy rosną do nieskończoności z prawdopodobieństwem $1$.
Tak też będziemy zakładać w całej tej pracy.
Dodatkowo będziemy zakładać, że firma wypłaca tzw. dywidendy.
Formalnie odpowiada to wprowadzeniu regulowanego procesu
\begin{equation}\label{eq:U} U^\pi_t = X_t - D^\pi_t, \q t \ge 0,
\end{equation}  gdzie  $\pi$
jest {\it dopuszczalną} strategią wypłaty (sposobem wypłaty lub też sposobem wyboru procesu $D^\pi_t$); patrz
de Finetti \cite{Fin}.
Cała ewolucja procesu w czasie jest możliwa do momentu ruiny
$$\tau^\pi =\inf\{t\in\mbb R_+: U^\pi_t < 0\}$$
czyli do pierwszego momentu, kiedy regulowany proces (po wypłacie dywidend) spadnie poniżej $0$.
Powyżej wspomniana {\it dopuszczalność} oznacza, że proces $D^\pi_t$ jest adaptowalny względem $\mathcal{F}_t$
(czyli decyzje są podejmowane w oparciu o dostępną wiedzę na temat dotychczasowej ewolucji procesu),
typu c\'adl\'ag oraz taki, który nie dopuszcza ruiny poprzez wypłatę dywidend:
\begin{equation}\label{exoruinliq} \text{dla każdego $t\leq\tau^\pi$ mamy}\q
\begin{cases}
D^\pi_{t}-D^\pi_{t-} \leq \le(X_t -
D^{\pi}_{t-}\ri)\vee 0, &  \\
\vspace{-0.5cm}\\
D^{\pi (c)}_t - D^{\pi (c)}_u \leq \E X_1 (t-u)
&
\text{$\forall u\in[0,t)$,}
\end{cases}
\end{equation}
gdzie $D^{\pi (c)}_t$ oznacza część ciągłą procesu $D^\pi_t$.
Oznaczmy przez $\Pi$ zbiór wszystkich dopuszczalnych strategii.

Naszym celem jest maksymalizacja wypłaty z punktu widzenia akcjonariuszy (dla firmy wyznacza ona zatem największą wypłatę jaką powinna wypłacić).
Jest ona równa średniej, zdyskontowanej i łącznej wypłacie dywidend wypłaconych do czasu ruiny:
\begin{equation}\label{value}
 v(x):=\sup_{\pi \in \Pi} \E_x\le[
\I_0^{\tau^\pi}\te{-qt}\td D^\pi_t\ri]. \end{equation}
Stała $q\geq 0$ jest stopą dyskontową.
Innymi słowy, naszym celem jest znalezienie $v$ oraz identyfikacja optymalnej strategii $\pi^*$ (o ile ona istnieje), która realizuje powyższe supremum.

\section{Uogólnienia}

Większość prezentowanego tutaj materiału po pewnych modyfikacjach może być użyta w ogólniejszych lub nieco innych
problemach dywidendowych. W szczególności można rozważać następujące wariacje prezentowanego modelu.

{\bf Koszty transakcji czyli tzw. kontrola impulsowa.}
Można wprowadzić tzw. stałe koszty transakcji $K>0$, które są płacone firmie za każdą wypłatę dywidend. W tym przypadku:
$$v(x) = \sup_{\pi\in\Pi} \E_x
\le[\I_0^{\tau^\pi}\te{-qt}\td D^\pi_t - K
\int_0^{\tau^\pi}\te{-qt}\td N^\pi_t\ri],
$$
gdzie $N^\pi_t$ jest procesem liczącym skoki procesu $D^\pi$. Mówi się wtedy o kontroli impulsowej związanej z wyborem momentów płacenia dywidend oraz
ich wielkości; patrz \cite{APP2015}.

{\bf Model z funkcją Gerbera-Shiu.}
W modelu tym uwzględnia się koszt wielkości długu w momencie ruiny. W tym celu wprowadza się tzw. {\it funkcję kary}
$w: (-\infty,0]\to (-\infty,0]$ oraz
funkcję Gerbera-Shiu:
$${\mathcal W}_w^\pi(x) := \E_x\le[\te{-q\tau^\pi}
w\le(U^\pi_{\tau^\pi}\ri); \tau^\pi<\infty \ri].$$
W przypadku kiedy $w = -1$ funkcja $-{\mathcal W}_w^\pi(x)$ opisuje tzw. prawdopodobieństwo ruiny.
W tym modelu funkcja wartości jest sumą:
$$ v(x):=\sup_{\pi \in \Pi} \left\{\E_x\le[
\I_0^{\tau^\pi}\te{-qt}\td D^\pi_t\ri]+\mathcal{W}(x)\right\};$$
patrz \cite{APP2015}.

{\bf Model dualny.} W tym przypadku rozważa się proces L\'evy'ego $X_t$ spektralnie dodatni, to znaczy pozbawiony ujemnych skoków; patrz \cite{dual} i referencje tam zawarte.

{\bf Inne momenty zatrzymania i wpłaty.} Można rozważać inne momenty ruiny. Szczególnie interesująca jest tzw. Paryska ruina, która dopuszcza przebywanie
poniżej $0$ procesu regulowanego ale nie dłużej niż ustalony z góry horyzont czasowy; patrz \cite{Irmina}. Są też modele, w których dopuszcza się wpłaty $R^\pi_t$ czyli dodatkową
regulację związaną z {\it 'podnoszeniem'} procesu w górę. Pierwszą możliwością są wpłaty kompensujące ruinę i utrzymujące proces w półosi $[0,\infty)$.
W tym przypadku proces regulowany jest równy
$$
U^{\pi}_t=X_t - D^{\pi}_t + R^{\pi}_t
$$
oraz funkcja wartości jest równa:
$$
v(x) = \sup_{\pi \in \Pi} \E_x\le[\I_0^\i \te{-qt}\td D^{\pi}_t -
\phi \I_0^\i \te{-qt}\td R^{\pi}_t\ri],
$$
gdzie $\phi$ są stałymi kosztami wpłat; patrz \cite{APP}.
Inną możliwością jest całkowanie powyżej tylko do momentu pojawienia się pierwszego skoku dostatecznie dużego aby regulowany proces przed wpłatą
znalazł się poniżej ustalonego poziomu $c<0$.

{\bf Proces odnowy skoków.}
W tym przypadku proces $N_t$ z (\ref{cramermod}) można zastąpić procesem odnowy, w którym wykładnicze
czasy pomiędzy kolejnymi skokami sa zastąpione dowolną, dodatnią zmienną losową. Literatura w tym przypadku jest ogromna, patrz np. \cite{AT09} i referencje tam zawarte.

{\bf $X$ jest rozwiązaniem stochastycznego równania różniczkowego.}
Najprostszym przypadkiem jest kiedy intensywność wpłat zależy od obecnego stanu rezerw, czyli $X_t$ rozwiązuje następujące stochastyczne równanie
różniczkowe:
$$dX_t=\eta(X_t)\;dt- dZ_t, \qquad X_0=x,$$
gdzie $Z_t$ jest procesem L\'evy'ego (najczęściej subordynatorem, czyli rosnącym procesem L\'evy'ego) zaś $\eta(\cdot)$ jest ustaloną, dodatnią funkcją, patrz \cite{ewa}.
W przypadku kiedy $Z_t=\sigma B_t$ dla $\sigma >0$, proces $X_t$ jest szczególnym procesem dyfuzji. Na ten temat jest bardzo aktywnie rozwijany i istnieje ogromna literatura, patrz np.
\cite{survey2} i referencje tam zawarte.

{\bf Inne funkcje wypłaty.}
Można rozważać tylko strategie $D^\pi_t$, które są absolutnie ciągłe z gęstością $d_t^\pi$. W tym przypadku wybiera się tzw. dodatnią funkcję użyteczności
$U(\cdot)$, która znika na ujemnej półosi. W tym przypadku:
$$
v(x) = \sup_{\pi \in \Pi} \E_x\le[\I_0^\i\te{-qt}U(d^{\pi}_t)\td t\ri];$$
patrz np. \cite{GHS}.

{\bf Wielowymiarowe procesy L\'evy'ego.} W tym przypadku rozważa się $X_t$ w $\mathbb{R}^d$ dla $d\geq 2$; patrz \cite{CP, AMP}.

Ze względu na przejrzystość dalsza część pracy będzie dotyczyć wyłącznie modelu zaprezentowanego w Rozdziale \ref{sec:model}.

\section{Dualna reprezentacja}
W przypadku kiedy $X_t$ jest procesem dyfuzji czasami udaje się pokazać dostateczną gładkość funkcji wartości $v$ tak aby móc stosować klasyczne metody
równań Hamiltona-Jacobiego-Bellmana (oznaczane dalej jako HJB).
%, do których wrócimy w Sekcji \ref{HJB}.
W przypadku procesów (skokowych) L\'evy'ego zwykle to się nie udaje, ponieważ
funkcja wartości jest ciągła (ale nie klasy $\mathcal{C}^1$) w przypadku kiedy $X_t$ jest procesem o ograniczonym wahaniu
oraz jest klasy $\mathcal{C}^1$ (ale nie klasy $\mathcal{C}^2$) kiedy proces $X_t$ ma nieograniczone wahanie.
Tymczasem równanie HJB jest następującej postaci:
\begin{eqnarray}\label{HJB}
\max\le\{{\mc L} g(x) - q g(x), 1 - g^\prime(x)\ri\} = 0, \qquad x>0,
\end{eqnarray}
gdzie %$g_+^\prime(x)$ jest prawostronną pochodną, oraz
$\mc L$ oznacza infitezymalny generator półgrupy przejścia procesu
$X$:
%\vspace{0.02cm}
\begin{equation}\label{eq:infgen}
\mc L g(x) = \frac{\s^2}{2}g''(x) + \eta g'(x) + \int_{\mbb R_+\backslash\{0\}}[g(x-y) - g(x) + yg'(x)\mathbf{1}_{\{y<1\}}]\nu(\td y),
\end{equation}
który działa na dostatecznie gładkich funkcjach $g$; patrz Sato~\cite[Tw. 31.5]{Sato}.
Powyżej miara $\nu$ jest miarą wartości bezwzględnej skoków (stąd jej nośnik jest na dodatniej półosi).
Można zaproponować dwie metody radzenia sobie z problemem niedostatecznej gładkości funkcji wartości:
metoda martyngałowa oraz metoda rozwiązań lepkościowych systemu HJB.
W tej pracy skoncentruję się na pierwszej metodzie.

Przypomnijmy, że proces $\{M_t, t\geq 0\}$ jest supermartyngałem (martyngałem) (względem filtracji $\mathcal{F}_t$)
jeśli $\E |M_t|<\infty$ oraz $\E[M_t|\mathcal{F}_s]\geq M_s$ ($\E[M_t|\mathcal{F}_s]=M_s$). Warunkowa wartość oczekiwana $\E[M_t|\mathcal{F}_s]$ to taka $\mathcal{F}_s$ całkowalna zmienna losowa $Z$, że
dla każdego zdarzenia $A\in\mathcal{F}_s$ mamy $\int_A Z\td \P=\int_A M_t\td \P$.
Niech $\tau_0^-=\inf\{t\geq 0: X_t<0\}$.

\begin{Def}\label{def:sss}\rm
Powiemy, że funkcja $g:\mathbb R_+\to\mathbb R$ jest {\it stochastycznym
nad-rozwiązaniem} równania HJB ~\eqref{HJB}
jeśli $g$ jest prawostronnie różniczkowalna z $g_+^\prime(x)\geq 1$ oraz proces
\begin{eqnarray}\label{smart}
M^g:=\left\{\te{-q\le(t\wedge \tau_0^-\ri)}g\le(X_{t\wedge \tau_0^-}\ri),\
t\geq 0\right\}\end{eqnarray}
jest jednostajnie całkowalnym nadmartyngałem ($M^g_t$ tworzą rodzinę jednostajnie całkowalnym zmiennych losowych).
Oznaczmy przez $\mc G$ rodzinę stochastycznych nad-rozwiązań.
\end{Def}

Okazuje się, że funkcja wartości $v$ jest najmniejszym stochastycznym superrozwiązaniem
(\ref{HJB}).
\begin{Thm}\label{thm:repg}
\begin{eqnarray}\label{eq:rep-g}
v(x) = \min_{g\in\mathcal G} g(x).
\end{eqnarray}
\end{Thm}

Zauważmy, że teraz w celu identyfikacji funkcji wartości formalnie wymagamy tylko istnienia
prawostronnej pochodnej i sprawdzenia warunku supermartyngałowego.
Konstrukcja ta przypomina konstrukcję tzw. obwiedni Snella znanej w matematyce finansowej.
Najważniejsze jest jednak poniższe twierdzenie weryfikacyjne będące oczywistym wnioskiem z powyższego twierdzenia.

\begin{Cor}\label{weryfikacja}
Jeśli dla pewnego $\pi^*$ funkcja $v_{\pi^*}(x)=\E_x\le[
\I_0^{\tau^{\pi^*}}\te{-qt}\td D^{\pi^*}_t\ri]$ jest stochastycznym nad-rozwiązaniem HJB (\ref{HJB}), to $v(x)=v_{\pi^*}(x)$ i $\pi^*$ jest optymalną strategią.
\end{Cor}

\section{Strategia barierowa}
Zgodnie z twierdzeniem weryfikacyjnym podanym we Wniosku \ref{weryfikacja}, zastosujemy tzw. metodę: "zgadnij i sprawdź", która polega
na znalezieniu funkcji wartości dla konkretnej strategii i sprawdzeniu czy jest ona stochastycznym superrozwiązaniem.
Tą wybraną strategią będzie strategia barierowa $\pi_a$, zgodnie z którą wypłacamy wszystko ponad ustalony poziom $a\geq 0$ jako dywidendy. Innymi słowy:
\begin{equation}\label{barrier}
D^{\pi_a}_t=\left(\overline{X}_t -a \right)\wedge 0,
\end{equation}
gdzie $\overline{X}_t:=\sup_{0\leq t}X_s$.
Niech
$$v_a(x):=\E_x\le[
I_q\ri],$$
gdzie  zmienna losowa
$$I_q:=\I_0^{\tau^{\pi_a}}\te{-qt}\td D^{\pi_a}_t$$
opisuje łączne zdyskontowane dywidendy wypłacone do momentu ruiny przy strategii $\pi_a$.
W następnym kroku zidentyfikujemy $v-a(x)$ tak aby znaleźć warunki dostateczne na to aby $v_a$ było stochastycznym superrozwiązaniem.

Zauważmy, że $v_a(x)=x-a+x_a(a)$ dla $x>a$, ponieważ kiedy $x>a$ od razu ( w  momencie $t=0$)
wypłacamy 'nadwyżkę' $x-a$ jak dywidendę i startujemy proces z $a$. Dodatkowo, dla $x\leq a$ aby móc zacząć wypłacać dywidendy
najpierw musimy się znaleźć w punkcie $a$, który przekraczamy w sposób ciągły (ponieważ nie mamy ściśle dodatnich skoków), unikając ujemnych wartości, gdzie następuje ruina.
Do tego czasu nie wypłacamy żadnych dywidend.
Stąd z własności Markowa dla procesu L\'evy'ego $X$ zachodzi
\begin{equation}
v_a(x)=h(x)v_a(a),
\end{equation}
gdzie
$$h_a(x):=\E_x\left[e^{-q\tau_a^+};\tau^+_a<\tau_0^-\right]$$
dla $\tau_a^+=\inf\{t\geq 0: X_t\geq a\}$.
Dość łatwo jest pokazać, że $h$ ma jednostronne pochodne oraz, że
$v_{a,-}^\prime(a)=v_{a,+}^\prime(a)=1$ co daje
$v_a(x)=h_a(x)/h_{a,-}^\prime(a)$.
Wystarczy zatem skoncentrować się na znalezieniu funkcji $h_a(x)$, która podaje rozwiązanie tzw. dwustronnego problemu wyjścia 'w górę'.
Ponownie z własności Markowa oraz braku dodatnich skoków wynika, że dla $b>a$ mamy $h_b(x)=h_a(x)h_b(a)$ i stąd istnieje funkcja $W^{(q)}(x)$ taka, że
\begin{equation}\label{exitup}
h_a(x)=\frac{W^{(q)}(x)}{W^{(q)}(a)}.\end{equation}
Funkcja $W^{(q)}(x)$ jest tzw. (pierwszą) funkcją skalująca i udaje się ją jednoznacznie zidentyfikować
za pomocą teorii martyngałów; patrz Dodatek \ref{Dodatek}.

W tym celu na mocy twierdzenia L\'evy'ego-Chińczyna dla procesu $X$ definiujemy wykładnik
Laplace'a $\psi(\theta$:
\begin{eqnarray}
 \label{eq:psi}\E[\te{\th( X_t-X_0)}] =
\te{t\psi(\th)}, \; \psi(\th)= \frac{\s^2}{2} \th^2 + \eta\,\th +
\int_0^\infty (
 \te{- \th x}- 1 + \theta x\mathbf{1}_{\{|x|<1\}}) \nu(\td x),
\end{eqnarray}
gdzie $\int_0^1 x^2\nu(\td x)<\infty$.
Wykładnik Laplace'a $\psi$ ze względu na brak skoków dodatnich procesu $X$
jest poprawnie zdefiniowany przynajmniej na dodatniej półosi,
na której jest ściśle wypukły spełniający
$\lim_{\th\to\infty}\ps(\th)=+\infty$. Dodatkowo, $\psi$ jest ściśle rosnąca na
 $[\F(0),\infty)$, gdzie $\F(0)$ jest największym pierwiastkiem równania
$\ps(\th)=0$. Przez
$\F:[0,\infty)\to[\F(0),\infty)$ będziemy oznaczać prawostronną funkcję odwrotną do $\psi$.
Wtedy funkcja
$W^{(q)}$ jest niemalejącą, ciągłą i znikającą na ujemnej półosi funkcją
posiadającą następującą transformatę Laplace'a:
\begin{eqnarray}
\label{Wq} \int_0^\infty \te{-\th x} W^{(q)}
(y) \td y = (\ps(\th) - q)^{-1},\q\q\th
> \F(q).
\end{eqnarray}
Powyższa funkcja pełni kluczową rolę w teorii fluktuacji.
Można pokazać, że jeśli
$X$ ma trajektorie o ograniczonym wahaniu wtedy $W^{(q)}(x)$ jest klasy $\mathcal{C}^1$
wtedy i tylko wtedy gdy miara $\nu$ nie ma atomów; patrz \cite{Hubalak}.
W przypadku kiedy $X$ ma trajektorie o nieograniczonym wahaniu
wtedy $W^{(q)}(x)$ jest zawsze klasy $\mathcal{C}^1$. Dodatkowo, jeśli $X$ ma komponentę brownowską, tzn. $\sigma >0$ to
$\mathcal{C}^1$ może być zastąpione przez $\mathcal{C}^2$.
%Dodatkowo, Song [13] udowodnił, że jeśli $v[x,\infty)$ is completely monotoniczna, to
%$W^{(q)}(x)$ jest klasy $\mathcal{C}^\infty$.
Załóżmy od teraz, że miara skoków nie ma atomów. Wtedy $W^{(q)}(x)\in\mathcal{C}^1$.
Podsumujmy:
\begin{Thm}
\begin{eqnarray}
\label{vax} v_a(x) &=&
\begin{cases} \frac{W^{(q)}(x)}{W^{(q)\prime}(a)},
& 0\leq x \leq a,\label{vax}\\
&\\
x - a + \frac{W^{(q)}(a)}{W^{(q)\prime}(a)}, & x >  a.
\end{cases}
\end{eqnarray}
\end{Thm}
Zauważmy, że gładkość $v_a(x)$ jest taka sama jak funkcji skalującej $W^{(q)}(x)$.
Niech $a^*$ będzie optymalną barierą, to znaczy maksymalizująca powyższą funkcję, co jest równoważne temu, że
\begin{equation}\label{astar}
a^* = \inf\{a>0: W^{(q)\prime}(a) \leq W^{(q)\prime}(x)\ \
\text{dla każdego $x\geq 0$}\}.
\end{equation}

Czasami wygodniej jest spojrzeć szerzej na strategię barierową.
Proces $\{a-U_t^{\pi_a}, t\geq 0\}$
względem $\P_a$ ma taki sam rozkład jak odbity proces od obecnego supremum
$\{(x\wedge Y_t, t \geq 0\}$ względem $\P$ dla $Y_t=:\overline{X}_t - X_t$.
Wtedy wypłacone dywidendy $D^{\pi_a}_t$ liczone do czasu ruiny $\tau^{\pi_a}$
mają taki sam rozkład jak $\overline{X}_t$ względem $\P$ liczony do pierwszego momentu
przekroczenia $a$ przez $Y$, to znaczy do $\inf\{t\geq 0: Y_t>a\}$.
Jeśli zatem  $\{(t, \epsilon_t) : t\geq 0\}$
jest (markowanym) punktowym procesem Poissona wycieczek z zera procesu
$Y$ zaś $L_t^{-1}$ związanym z nim czasem lokalnym, to po zamianie czasu mamy:
 \[
I_q = \int_0^\infty e^{-qL^{-1}_t}\mathbf{1}_{(\sup_{s\leq
t}\overline{\epsilon}_s \leq a)}dt,
\]
gdzie $\overline{\epsilon}_t$ jest wysokością wycieczki. Stąd $I_q$ jest szczególnym przypadkiem następującej całki:
 \begin{equation}
I_q=\int_0^\infty e^{-q\xi_t}\mathbf{1}_{(\sup_{s\leq t}\Delta \beta_s
\leq a)}dt, \label{abstraction}
\end{equation}
gdzie $\Delta\beta_t = \beta_t-\beta_{t-}$ oraz
$(\xi,\beta)= \{(\xi_t , \beta_t) : t\geq 0\}$ jest (tutaj szczególnym) dwuwymiarowym subordynatorem.

Z powyższej obserwacji można uzyskać $\E_xI^n_q$ dla $n=1,2,3,\ldots$. Dodatkowo z
twierdzeń Tauberwoskich typu Kashary (patrz \cite[Tw. 4.12]{bingham} ) można uzyskać na przykład następujący rezultat.

\begin{Thm}\label{rivero}
Załóżmy, że $X$ spełnia warunek Spitzera-Doneya w $0$, tzn. istnieje
$\rho\in(0,1)$ takie, że
\begin{equation}\label{spitzer0}
\lim_{t\downarrow0}\frac{1}{t}\int_0^t \P(X_s \geq 0)ds =
\rho.
\end{equation}
Wtedy
\[
\lim_{y\to\infty}\frac{-\log \mathbb{P}_x(I_q >y)}{(1-\rho)\varpi(y)}=1,
\]
gdzie  $\varpi(x)$ jest jedynym rozwiązaniem w $(0,\infty)$ równania $\psi(\theta)=\theta x$.
\end{Thm}

Warunek (\ref{spitzer0}) jest równoważny następującemu warunkowi (patrz \cite{Doney}):
\begin{equation}\label{spitzer}
\lim_{t\downarrow0} \mathbb{P}(X_t \geq 0)=\rho.
\end{equation}
Procesem spełniającym powyższy warunek jest na przykład spektralnie ujemny proces stabilny z indeksem
$\alpha\in(1,2)$.

\section{Optymalność i rozwiązania lepkościowe}

Zauważmy, że wprost z definicji $a^*$ wynika, że
$v_{a^*}^\prime(x)\geq 1$. Dodatkowo, z (\ref{exitup}) wynika, że
$$\E\left[e^{-q \tau_0^-\wedge \tau_a^+}W^{(q)}(X_{\tau_0^-\wedge \tau_a^+})|\mathcal{F}_t\right]
=\te{-q(t\wedge \tau_0^-\wedge \tau_a^+)}W^{(q)}(X_{t\wedge \tau_0^-\wedge \tau_a^+})$$
jest jednostajnie całkowalnym martyngałem.
Stąd proces $\te{-qt\wedge
\tau_0^-\wedge \tau_a^+}v_{a^*}(X_{t\wedge \tau_0^-\wedge \tau_a^+})$
dziedziczy własność martyngałową.
Wreszcie dla $x>a^*$ z (\ref{eq:infgen}) mamy, że
\begin{equation}\label{wargen}
\mc L v_{a^*}(x) = \eta + \int_{\mbb R_+\backslash\{0\}}[v_{a^*}(x)(x-y) - v_{a^*}(x) + y\mathbf{1}_{\{y<1\}}]\nu(\td y)
\end{equation}
niezależnie od gładkości funkcji skalującej.
Z twierdzenia Dynkina zastosowanego do wykładniczo zabijanego procesu $X$ z intensywnością $q$
oraz z Wniosku \ref{weryfikacja} wynika następujący rezultat.

\begin{Thm}
Jeśli
\begin{equation}\label{mainzal1}
\mc L v_{a^*}(x)-q v_{a^*}(x)\leq 0 \qquad \text{dla $x>a^*$,}
\end{equation}
to $v(x)=v_{a^*}(x)$ oraz $\pi_{a^*}$ jest optymalną strategią $\pi^*$.
\end{Thm}
\begin{Rem}\label{Uwaga1}\rm
Z powyższych rozważań wynika również, ze strategia barierowa $\pi_{a^*}$ jest zawsze optymalna w podrodzinie strategii ograniczonych przez $a^*$.
\end{Rem}

Zanim przyjrzymy się założeniu (\ref{mainzal1}),
zastanówmy się jaka strategia jest zawsze optymalna (bez dodatkowych założeń).
Okazuje się, że tą strategią jest tzw. strategia bandowa zdefiniowana w następujący sposób.
Strategia $\pi_{\unl b,\unl a}$ jest związana z ciągiem
$$ b_1=0 \leq a_1 < b_2 \leq a_2< ...
<b_{n} \leq a_n < \ldots,
$$ i następującą regulacją:
\begin{itemize}
\item[(i)] Kiedy $U^{\pi_{\underline{b},\underline{a}}}_t=U^{\underline{b},\underline{a}}_t=
y \in (a_n, b_{n+1})$,
wypłacamy $y-a_n$;

\item[(ii)] Kiedy $U^{\underline{b},\underline{a}}=a_n$
odbijamy proces $X$ od bariery $a_n$;

\item[(iii)] Kiedy $U^{\underline{a},\underline{b}} \in [b_n,a_{n})$
nie wypłacamy dywidend.
\end{itemize}

Mamy więc trzy obszary:
$\mathcal{A}_1$ gdzie wypłacamy dywideny w postaci jednorazowej wypłaty (co odpowiada skokowi procesu regulowanego do najbliższej 'odbijającej' bariery $a_n$),
$\mathcal{A}_2$ gdzie proces regulowany jest odbijany w $a_n$ i $\mathcal{A}_3$ gdzie nie są wypłacane dywidendy (brak regulacji).
Przypadek $b_1=0$ i $a_1=a^*$ odpowiada strategii barierowej.
Argumenty oparte o indukcję matematyczną oraz obserwacji podanej w Uwadze \ref{Uwaga1} dają (po dość zaawansowanej analizie) następujące twierdzenie.

\begin{Thm}
Strategia bandowa jest zawsze optymalna.
\end{Thm}

Z tą strategią nieodłącznie związane jest rozwiązanie lepkościowe systemu
(\ref{HJB}). W przypadku kiedy można pokazać, że funkcja wartości $v$ jest wystarczająco gładka
to równanie staje się zwykłym równaniem różniczkowo-całkowym.
Zacznijmy od tego, że dość łatwo jest pokazać, że funkcja
$v(x)$ jest zawsze lokalnie Lipschitzowska.

\begin{Def}
\label{Viscosity}\rm Funkcję lokalnie Lipschitzowską $\overline{u}: \mathbb{R}%
_{+}\rightarrow\mathbb{R}$\ nazywamy \textit{nad-rozwiązaniem lepkościowym}
(\ref{HJB})\ w $x \in\mathbb{R}_{+}$\ jeśli każda funkcja klasy $\mathcal{C}^2$
$\varphi:\mathbb{R}_{+}\rightarrow\mathbb{R}%
\ $ taka, że $\varphi(x)=\overline{u}(x)$ oraz
$\overline{u}-\varphi$\ osiąga minimum w $x$, spełnia następującą nierówność:
\[
\max\le\{{\mc L} \varphi(x) - q \varphi(x), 1 - \varphi^\prime(x)\ri\} \leq 0.
\]
Podobnie, funkcję lokalnie Lipschitzowską $\underline{u}: \mathbb{R}%
_{+}\rightarrow\mathbb{R}$\ nazywamy \textit{pod-rozwiązaniem lepkościowym}
(\ref{HJB})\ w $x \in\mathbb{R}_{+}$\ jeśli każda funkcja klasy $\mathcal{C}^2$
$\vartheta:\mathbb{R}_{+}\rightarrow\mathbb{R}%
\ $ taka, że $\vartheta(x)=\underline{u}(x)$ oraz
$\underline{u}-\varphi$\ osiąga minimum w $x$, spełnia następującą nierówność:
\[
\max\le\{{\mc L} \vartheta(x) - q \vartheta(x), 1 - \vartheta^\prime(x)\ri\} \geq 0.
\]
Funkcja, która jest jednocześnie pod-rozwiązaniem i nad-rozwiązaniem lepkościowym będziemy nazywać rozwiązaniem lepkościowym
(\ref{HJB}).
\end{Def}

\begin{Thm}
\label{Prop V is a viscosity supersolution} Funkcja $v$ jest rozwiązaniem lepkościowym
systemu (\ref{HJB}).
\end{Thm}

Mając to rozwiązanie można również zidentyfikować bandową strategię optymalną:
\begin{eqnarray*}
\mathcal{A}_1&=&\left\{x\in\mathbb{R}_{+}:\mathcal{L}%
v(x) -qv(x) <0\text{, }1-v^\prime(x)=0\right\}\\
\mathcal{A}_2&=&\left\{x\in\mathbb{R}_{+}:\mathcal{L}%
v(x) -qv(x) =0\right\}\\
\mathcal{A}_3&=&\left(\mathcal{A}_1\cup \mathcal{A}_3\right)^c.
\end{eqnarray*}

Powyższy rezultat ma ogromną wagę choćby ze względu na analizę numeryczną.
Jednakże nawet jeśli wiemy jaka strategia jest zawsze optymalna
intrygujące pozostaje pytanie kiedy najprostsza
strategia barierowa $\pi_{a^*}$ jest optymalna, lub, innymi słowy, kiedy zachodzi warunek
(\ref{mainzal1}). Oczywiście optymalne warunki powinny być wyrażone poprzez wyjściowe parametry modelu $(\eta, \sigma, \mu)$ czyli przez tzw. {\it trójkę L\'evy'ego}.
Korzystając z postaci funkcji $v_{a^*}(x)$ podanej w (\ref{vax}) oraz
(\ref{wargen}) można znaleźć poprzez dość proste rachunki zaskakujący warunek dostateczny
na to aby zachodziła nierówność  (\ref{mainzal1}).

\begin{Thm}
Jeśli
\begin{equation}\label{suff1}
W^{(q),\prime}(x)\geq W^{(q),\prime} (y)\qquad \text{dla każdego $x>y>a^*$,}
\end{equation}
to $v(x)=v_{a^*}(x)$ oraz $\pi_{a^*}$ jest optymalną strategią.
\end{Thm}
Zgodnie z definicją $a^*$ podaną w (\ref{astar}) dostateczny warunek (\ref{suff1}) mówi, że strategia barierowa jest optymalna jeśli
funkcja skalująca jest wypukła
po osiągnięciu globalnego minimum swojej pochodnej.
Ta własność funkcji skalującej jest związana z innym ciekawym fenomenem opartym o teorię potencjału subordynatorów, który opiszemy w następnym rozdziale.

\section{Funkcja skalująca i funkcje Bernsteina}
Zacznijmy od tego, że $q$-ta funkcja skalująca wyraża się poprzez zerową funkcję
skalującą na nowej mierze probabilistycznej $\P^{\Phi(q)}$ zdefiniowanej poprzez pochodną Radona-Nikodyma typu
Girsanowa:
\begin{equation}\label{girsanov}
	\frac{\td \P_{x}^{\Phi(q)}}{\td \P_{x}}\biggr|_{\mathcal{F}_{t}}
	=\mathcal{E}_t(\Phi(q)),
\end{equation}
gdzie
\begin{equation}\label{wyklmart}\mathcal{E}_t(\theta):=e^{\theta (X_t-X_0)-\psi(\theta)t}\end{equation}
jest wykładniczym martyngałem.
Przypomnijmy, że $\Phi(q)=\psi^{-1}(q)$ dla wykładnika Laplace'a $\psi$ procesu $X$ oraz ze względu na założenie, że proces $X$ zmierza do nieskończoności, mamy
$\Phi(0)>0$.
Na $\P^{\Phi(q)}$ proces $X$ pozostaje spektralnie ujemnym procesem L\'evy'ego z nowym wykładnikiem Laplace'a
\begin{equation}\label{newpsi}
\psi_q(\theta)=\Psi(\theta+\Phi(q))-q.\end{equation}
Formalnie mamy:
\begin{equation}\label{scalerepr1}
W^{(q)}(x)=e^{\Phi(q) x}W_{\Phi(q)}(x),
\end{equation}
gdzie $W_{\Phi(q)}(x)=W_{\Phi(q)}^{(0)}(x)$ jest zerową funkcją skalującą na $\P^{\Phi(q)}$ co z (\ref{Wq}) jest równoważne
faktowi, że
\begin{equation}\label{zero}
\int_0^\infty \te{-\th x} W_{\Phi(q)}
(y) \td y = (\psi_q(\th))^{-1};
\end{equation}
patrz \cite{KypPal}.
Z procesem odbitym w infimum $\widehat{Y}_t=X_t-\inf_{s\leq t}X_s$ można związać czas lokalny w $0$, $\{\widehat{L}_t, t\geq 0\}$
oraz rozważać schodzący proces drabinowy $\{(\widehat{L}_t^{-1}, \widehat{H}_t), t\geq 0\}$, gdzie $\widehat{H}_t:=X_{\widehat{L}^{-1}_t}$ dla $t<\widehat{L}_\infty$.
Dla tego (zabijanego wykładniczo) dwuwymiarowego subordynatora możemy zdefiniować jego wykładnik Laplace'a $\widehat{\kappa}$ w następujący sposób:
$$
\E\left[e^{-\theta\widehat{L}_1^{-1}-\beta \widehat{H}_1}\right]=e^{-\kappa(\theta, \beta)}.
$$
Ze względu na brak dodatnich skoków wchodzący proces drabinowy $H_t:=X_{L^{-1}_t}$ dla czasu lokalnego w supremum $L_t$ jest liniowy.
Z faktoryzacji Wienera-Hopfa (patrz \cite{Sato} i \cite{Kyprianou}) na $\P^{\Phi(q)}$ mamy:
\begin{equation}\label{WH}
\psi_q(\theta)=\theta\widehat{\kappa}(q,\Phi(q)+\theta).
\end{equation}
Stąd:
\begin{equation}\label{scalerepr0}
W_{\Phi(q)}(x)=\E^{\Phi(q)}\int_0^\infty \mathbf{1}_{\{\widehat{H}_t\in [0,x]\}}\td t
\end{equation}
jest miarą potencjałową (funkcją odnowy) subordynatora $\widehat{H}$ z wykładnikiem Laplace'a $\widehat{\kappa}(q,\cdot)$ na $\P^{\Phi(q)}$.

Załóżmy teraz, że
na $\P$ miara skoków $\nu$ ma gęstość $\varrho(x)$ (tzn.
$\nu([0,y])=\int_0^y \varrho(z)\td z$),
która jest log-wypukła.

Wtedy funkcja $x\rightarrow \nu(x,\infty)$ jest również log-wypukła.
Z \cite[Lem. 2.3]{krs} miara skoków drabinowego $\widehat{H}$ na $\P^{\Phi(q)}$ ma miarę $\Upsilon$ równą
\begin{equation}\label{ladderjump}
\Upsilon(x,\infty):=e^{\Phi(q)x}\int_x^\infty e^{-\Phi(q)y}\nu(y,\infty)\td y.
\end{equation}
Funkcja $x\rightarrow \Upsilon(x,\infty)$ jest log-wypukła.
Z \cite[Tw. 2.4]{Song} i \cite[Lem. 2.2]{krs} wynika, że
$\widehat{\kappa}(q,\cdot)$ jest funkcją Bernsteina  z nierosnącą i {\it wypukłą} potencjałową gęstością $u_q(x)$
spełniającą równanie Volterry:
$$d^* u_q(x) +\int_0^t \Upsilon(x-y, \infty)u_q(y)\td y=1, \qquad x >0,$$
gdzie $d^*$ jest dryfem $\widehat{H}_t$, tzn. $d^*=\lim_{\theta\to\infty} \widehat{\kappa}(q,\theta)/\theta$.
Zatem  $\widehat{\kappa}(q,\Phi(q)+\cdot)$ jest również funkcją Bernsteina  z gęstością $e^{-\Phi(q)x}u_q(x)$, to znaczy, że
$$\frac{1}{\widehat{\kappa}(q,\Phi(q)+\theta)}=\int_0^\infty e^{-\theta x} e^{-\Phi(q)x}u_q(x)\td x.$$
Na mocy wzorów (\ref{zero}), (\ref{WH}) oraz (\ref{scalerepr0}) mamy, że
\begin{equation}\label{scalerepr0b}
W_{\Phi(q)}(x)=d +\int_0^x e^{-\Phi(q)y}u_q(y)\td y,
\end{equation}
gdzie $d=\lim_{\theta \to \infty} 1/\widehat{\kappa}(q,\Phi(q)+\theta)\geq 0$.

Korzystając z (\ref{scalerepr1}) oraz  (\ref{scalerepr0b}) można udowodnić następujące twierdzenie.
\begin{Thm}\label{optimalmain}
Załóżmy, że gęstość miary skoków procesu L\'evy'ego $X$ ma log-wypukłą gęstość.
Wtedy $v(x)=v_{a^*}(x)$ i strategia barierowa jest optymalna.
\end{Thm}

\begin{Rem}\rm
Ponieważ każda funkcja kompletnie monotoniczna jest log-wypukła, więc
jeśli miara skoków ma kompletnie monotoniczną gęstość, to strategia barierowa jest optymalna.
Założenie to pojawia się też w innych rezultatach związanych najczęściej z procesem supremum $\overline{X}_t$ lub
wyjściem procesu $X$ z pewnego obszaru (najczęściej półprostej); patrz
np. \cite{Kwasnicki} i inne prace tego autora.
\end{Rem}

\section{Przykłady}
{\bf Liniowy ruch Browna.}
Jeśli $X_t = \s B_t + \eta t$ dla $\sigma >0$, to
$$
W^{(q)}(x) = \frac{1}{\s^2\d} [\te{(-\omega+\d) x} -
\te{-(\omega+\d)x}],
$$
gdzie $\d = \s^{-2}\sqrt{\eta^2 + 2q\s^2}$ i $\omega = \m/\s^2$.
Wtedy z Tw. \ref{optimalmain} strategią optymalną jest strategia barierowa $\pi_{a^*}$ z
$$
a^* = \log \le|\frac{\d + \omega}{\d - \omega}\ri|^{1/\d}.
$$

{\bf Model Cram\'{e}ra-Lundberga z wykładniczymi skokami.}
Jeśli $X$ jest dany w (\ref{cramermod}), gdzie
$C_i$ mają rozkład wykładniczy z parametrem $\mu$, to
$\psi(\th)=(\eta+\E C_1)\th-\lambda \th/(\m+\th)$
i$$
W^{(q)}(x) = p^{-1}\left(A_+ \te{q^+(q)x} - A_-
\te{q^-(q)x}\right),
$$
gdzie $A_\pm = \frac{\mu + q^\pm(q)}{q^+(q)-q^-(q)}$ z
$q^+(q)=\F(q)$ i $q^-(q)$ będącymi pierwiastkami kwadratowego równania $\psi(\th)=q$.
Wtedy z Tw. \ref{optimalmain} strategią optymalną jest strategia barierowa $\pi_{a^*}$ z
$$a^* =
\frac{1}{q^+(q) - q^-(q)}
\log \frac{q^-(q)^2(\m+q^-(q))}{q^+(q)^2(\m+q^+(q))}.$$

{\bf Przykład Azcue-Muler podany w \cite{AM}.}
Jest to przykład procesu, dla którego strategia barierowa nie jest optymalna.
Jeśli $X$ jest dany w (\ref{cramermod}), gdzie
$C_i$ mają rozkład Gamma rzędu $2$ z parametrem
$\mu=1$, tzn.
$$\P(C_1\in \td y)=ye^{-y}\td y.$$
Czyli gęstość skoków nie jest log-wypukła.
Dodatkowo niech $\lambda=10$, $\eta =\frac{97}{5}$, $q =\frac{1}{10}$.
%$\th =\frac{7}{100}$.
Wtedy
$\psi (\theta) - q = (\eta+\E C_1)\theta  + \lambda (\frac{\mu}{\mu+\theta})^2 -\lambda-q
=\frac{\eta+\E C_1}{(\mu+\theta)^2}(\theta +\zeta_1)(\theta +\zeta_2)(\theta -\zeta_0),$
gdzie $\zeta_0 \approx 0.0396$, $\zeta_1\approx 0.0794$, $\zeta_2 \approx 1.4882$.
Stąd:
$$W^{(q)}(x)=\sum_{i=1}^3 e^{\zeta_i x}A_i$$
dla pewnych, jawnych stałych $A_i$. Optymalną strategią $\pi^*$ jest strategia bandowa z bandami $a_1
\approx 1.83$ i $b_1\approx 10.45$.
% czyli funkcja wartości równa się:

\section*{Podziękowania}
Praca powstała w wyniku realizacji projektu badawczego o nr 2013/09/B/HS4/01496
finansowanego ze środków Narodowego Centrum Nauki.

\section{Dodatek}\label{Dodatek}
W tym dodatku (w oparciu o \cite{KypPal})
krótko opiszemy jak za pomocą a teorii martyngałów rozwiązuje się
problem wyjścia (\ref{exitup}):
$$h_a(x)=\E_x\left[e^{-q\tau_a^+};\tau^+_a<\tau_0^-\right]
=\frac{W^{(q)}(x)}{W^{(q)}(a)}.$$

Główną rolę spełnia tzw. martyngał Kella--Whitta (lub martyngał Kennedy'ego):
\begin{equation}
	K_{t}:=\psi(\theta)\int_0^{t}
	\E^{-\theta Y_{s}}\,\td s
	+1-\E^{-\theta Y_{t}}-\theta \overline{X}_{t},
\label{kyprianou-palmowski:hides-WH}
\end{equation}
którego martyngałowość wynika z
równości:
\begin{eqnarray*}
	\td K_{t} &=&-e^{-\theta \overline{X}_{t}+\psi(\theta)t}
	\,\td \mathcal{E}_{t}(\theta)
\end{eqnarray*}
dla $\mathcal{E}_{t}(\theta)$ podanego w (\ref{wyklmart}).
Niech $\mathbf{e}_{q}$ będzie niezależną od procesu $X$ zmienną losową z parametrem $q$.
Zauważmy,że stosując twierdzenie Dooba o optymalnym momencie zatrzymania
mamy:
\begin{equation*}
	\E_{x}\biggl(\frac{\mathcal{E}_{t\wedge \tau_a^+}(\Phi(q))}
	{\mathcal{E}_0(\Phi(q))}\biggr)
	=\E_{x}\Bigl(e
^{\Phi(q) (X_{t\wedge\tau_a^+}-x)
	-q(t\wedge \tau_a^+)}\Bigr) =1.
\end{equation*}
Czyli
\begin{equation}\label{kyprianou-palmowski:palmowski}
	\P_x(\overline{X}_{\mathbf{e}_{q}}>a)=\P(\tau_a^+<\mathbf{e}_{q})=\E_{x}\bigl(e^{-q\tau_a^+};\tau_a^+<\infty \bigr)
	=e^{-\Phi(q) (a-x)}.
\end{equation}
Stąd zmienna losowa $\overline{X}_{\mathbf{e}_{q}}$ na $\P$ ma rozkład wykładniczy z parametrem $\Phi(q)$.
Ważna jest też w geometryczna obserwacja, że
$-\underline{X}_t:=\inf_{s\leq t}X_s$ ma taki sam rozkład jak $Y_t$.
Ponownie stosując twierdzenie Dooba o optymalnym momencie zatrzymania
dla $K_t$ używając powyższej obserwacji uzyskujemy:
\begin{equation}
	\E\bigl(\E^{\theta\underline{X}_{\mathbf{e}_{q}}}\bigr)
	=\frac{q(\theta -\Phi(q))}{\Phi(q)(\psi(\theta) -q)}\,.
\label{kyprianou-palmowski:I-factor}
\end{equation}

Załóżmy na początku, że
$\psi'(0^+)>0$ czyli $X_t \to\infty$ prawie wszędzie. Wtedy zmienna losowa
$\ -\underline{X}_{\infty}$ jest skończona prawie wszędzie.
Wzięcie $q\downarrow 0$ w (\ref{kyprianou-palmowski:I-factor}) daje:
\begin{equation*}
	\E\bigl(e^{\theta \underline{X}_{\infty}}\bigr)
	=\psi'(0)\,\frac{\theta}{\psi(\theta)}.
\end{equation*}
Całkowanie przez części produkuje:
\begin{eqnarray*}
	\E\bigl(e^{\theta \underline{X}_{\infty}}\bigr)
	&=&\int_{[0,\infty)}
	e^{-\theta x}\,\P(-\underline{X}_{\infty}\in \td x)
\\
	&=&\theta\int_0^{\infty}e^{-\theta x}
	\P(-\underline{X}_{\infty}<x)\,\td x
\\
	&=&\theta \int_0^{\infty}e^{-\theta x}
	\P_{x}(\underline{X}_{\infty}>0)\,\td x.
\end{eqnarray*}
Definiując zerową funkcję skalującą w następujący sposób:
\begin{equation}
\label{kyprianou-palmowski:forces-left-cts}
	W^{(0)}(x)=W(x)=\frac{1}{\psi'(0^+)}\,P_{x}(\underline{X}_\infty > 0).
\end{equation}
łatwo sprawdzić, ze spełnia ona definicję: jest lewostronnie ciągła, znika na ujemnej półosi i
ma transformatę Laplace'a $1/\psi(\theta)$ dla wszystkich
$\theta>0$. Korzystając z mocnej własności Markowa oraz z faktu, że proces L\'evy'ego jest spektralnie ujemny (brak dodatnich skoków) uzyskujemy:
\begin{eqnarray}
	\hbox to 1em{$\P_{x}(\underline{X}_{\infty} > 0)$\hss}
	&&\nonumber\\
	&=&\E_{x}\,\P_{x}(\underline{X}_{\infty} >0\,|\,\mathcal{F}_{\tau_a^+})
	\nonumber\\
	&=&\E_{x}\Bigl(
	\P_a(\underline{X}_{\infty} > 0);\tau_a^+<\tau_0^-\Bigr)
	+\E_{x}\Bigl(
	\P_{X_{\tau_0^-}}(\underline{X}_{\infty}> 0);\tau_a^+>\tau_0^-)\Bigr)
\label{kyprianou-palmowski:importance of left cts}\\
	&=&\P_a(\underline{X}_{\infty} >0)\,\P_{x}(\tau_a^+<\tau_0^-).
	\nonumber
\end{eqnarray}
Czyli:
\begin{equation}
	\P_{x}(\tau_a^+<\tau_0^-)=\frac{W(x)}{W(a)}.
\label{kyprianou-palmowski:q=0}
\end{equation}

Załóżmy teraz, że $q>0$ lub $\psi'(0)<0$ i $q=0$
Z wypukłości
wykładnika Laplace'a $\psi$
wynika, że $\psi_{q}'\left(0\right) =\psi^{\prime}\left(\Phi
(q) \right) >0$ dla $\psi_q$ zdefiniowanego w (\ref{newpsi}).
Zamiana
miary (\ref{girsanov}) daje teraz ogólne rozwiązanie:
\begin{eqnarray*}
	\E_{x}\Bigl(e^{-q\tau_a^+};\tau_a^+<\tau_0^-\Bigr)
	&=&\E_{x}\biggl(\frac{\mathcal{E}_{\tau_a^+}(\Phi(q))}
	{\mathcal{E}_0(\Phi(q))}\mathbf{1}_{(\tau_a^+<\tau_0^-)}\biggr)
	e^{-\Phi(q)(a-x)}
\\
	&=&e^{-\Phi(q)(a-x)}\P_{x}^{\Phi(q)}
	\bigl(\tau_a^+<\tau_0^-\bigr)=\frac{W^{(q)}(x)}{W^{(q)}(a)}
\end{eqnarray*}
dla $q$-tej funkcji skalującej $W^{(q)}$ podanej w (\ref{scalerepr1}).
Co więcej, dla $\theta >\Phi(q)$ mamy:
\begin{eqnarray}
	\int_0^{\infty}e^{-\theta x}W^{(q)}(x)\,\td x
	&=&\int_0^{\infty}e^{-(\theta -\Phi(q))x}W_{\Phi(q)}(x)\,\td x
	\notag
	\\
	&=&\frac{1}{\psi_{\Phi(q)}(\theta -\Phi(q))} \notag
	\\
	&=&\frac{1}{\psi(\theta)-q}\,,
	\label{kyprianou-palmowski:WqLT}
\end{eqnarray}
co jest zgodne z definicją (\ref{Wq}).

\end{document}